\documentclass{amsart}

\usepackage{amsmath,amssymb,amsthm,pinlabel}

\def\co{\colon\thinspace}

\newcommand{\Diff}{\mbox{\rm Diff}}
\newcommand{\Cont}{\mbox{\rm Cont}}

\newcommand{\SO}{\mbox{\rm SO}}

\newcommand{\D}{\mathcal{D}}
\newcommand{\N}{\mathbb{N}}
\newcommand{\R}{\mathbb{R}}
\newcommand{\Z}{\mathbb{Z}}
\newcommand{\Q}{\mathbb{Q}}
\newcommand{\C}{\mathbb{C}}
\newcommand{\LL}{\mathbb{L}}
\newcommand{\calL}{\mathcal{L}}

\newcommand{\tb}{{\tt tb}}

\newcommand{\lk}{{\tt lk}}

\newcommand{\xist}{\xi_{\mathrm{st}}}

\newcommand{\rme}{{\mathrm{e}}}
\newcommand{\rmd}{{\mathrm{d}}}
\newcommand{\rmD}{{\mathrm{D}}}
\newcommand{\bfi}{{\mathbf{i}}}
\newcommand{\bfp}{{\mathbf{p}}}
\newcommand{\bfq}{{\mathbf{q}}}
\newcommand{\frakp}{{\mathfrak{p}}}

\newtheorem{thm}{Theorem}

\newtheorem{prop}[thm]{Proposition}
\newtheorem{cor}[thm]{Corollary}

\theoremstyle{definition}
\newtheorem*{rem}{Remark}

\newtheorem*{defn}{Definition}
\newtheorem*{ex}{Example}
\newtheorem*{exs}{Examples}

\begin{document}

\title{Contact structures and geometric topology}

\author{Hansj\"org Geiges}
\address{Mathematisches Institut, Universit\"at zu K\"oln,
Weyertal 86--90, 50931 K\"oln, Germany}
\email{geiges@math.uni-koeln.de}
\thanks{Some parts of the research described in this survey were
supported by DFG grant GE 1245/1 within the framework of the
Schwer\-punkt\-pro\-gramm 1154 ``Globale
Differentialgeometrie''.}

\date{}



\maketitle

\section{Introduction}
A {\em contact structure\/} on a manifold $M$ of dimension $2n+1$ is a
tangent hyperplane field, i.e.\ a $2n$-dimensional
sub-bundle $\xi$ of the tangent bundle~$TM$, satisfying the
following {\em maximal non-integrability\/} condition: if $\xi$
is written locally as the kernel of a differential $1$-form~$\alpha$,
then $\alpha\wedge (\rmd\alpha )^n$ is required to be nowhere zero on its
domain of definition. Notice that $\xi$ determines
$\alpha$ up to multiplication by a smooth nowhere zero function~$f$.
So the contact condition is independent
of the choice of $1$-form defining~$\xi$, since
$(f\alpha)\wedge \rmd (f\alpha )^n=f^{n+1}\alpha\wedge (\rmd\alpha )^n$.
I shall always assume our contact structures to be {\em coorientable},
which is equivalent to saying that we can write $\xi =\ker\alpha$
with a $1$-form $\alpha$ defined on all of~$M$; such an $\alpha$
is called a {\em contact form}. Then
$\alpha\wedge (\rmd\alpha )^n$ is a volume form on~$M$,
so a {\em contact manifold\/} $(M,\xi =\ker\alpha )$ has to be orientable.

The classical Darboux theorem states that any contact form $\alpha$
can locally be written, in suitable coordinates, as
$\alpha =\rmd z+\sum_{i=1}^nx_i\, \rmd y_i$. This is one of the reasons
why the most interesting aspects of contact geometry are of global nature.

Contact structures provide the mathematical language for many phenomena
in classical mechanics, geometric optics and thermodynamics.
Equally important for the interest in these structures are
their relations with symplectic, Riemannian and complex geometry.
These aspects are surveyed in~\cite{geig01} and~\cite[Chapter~1]{geig08}.

In the last two decades it has become increasingly apparent that
contact manifolds constitute a natural framework for many
problems in low-dimensional geometric topology. As hypersurfaces in
symplectic $4$-manifolds, $3$-dimensional contact manifolds
build a bridge to $4$-manifold topology. This interplay
between dimensions three and four has helped solve some
long-standing problems in knot theory. One salient example is
the result of Kronheimer--Mrowka that all non-trivial knots
in the $3$-sphere $S^3$ have the so-called property~P;
see Section~\ref{section:caps} below. Their proof is
based on a result of Eliashberg and, independently, Etnyre
that any symplectic filling of a
$3$-dimensional contact manifold can be capped off to a closed symplectic
$4$-manifold.

Moreover, contact topology has inspired new approaches to some known
results. Pride of place has to be given to Eliashberg's proof~\cite{elia92}
of Cerf's theorem that any diffeomorphism of $S^3$ extends to the
$4$-ball, based on the classification of contact structures
on~$S^3$; see~\cite{geig08} for an exposition of that proof.

Arguably the most influential contact topological result of the last
decade is due to Giroux~\cite{giro02}, cf.~\cite{etny06}
and~\cite{coli08}. He established a correspondence between
contact structures on a given manifold and open book decompositions
of that manifold; in dimension three and subject to suitable
equivalences on either set of structures, this correspondence is
actually one-to-one.

In the present article I want to
survey a selection of these recent developments in contact
topology. In Section~\ref{section:basic}
a few basic contact geometric concepts will be reviewed.
I then discuss some of the results highlighted above, and others
besides, from a somewhat idiosyncratic point of view.
As starting point I take a surgery presentation
of contact $3$-manifolds due to Fan Ding and yours truly;
this is the content of Section~\ref{section:surgery}.

In Section~\ref{section:applications} we then turn to
applications of this structure theorem.
For instance, one can use it
to derive an adapted open book decomposition
(see Section~\ref{section:surgery-book}),
thus providing an alternative proof for one direction of the
Giroux correspondence in dimension three.
In Section~\ref{section:caps}
I shall also explain in outline
how symplectic caps can be constructed directly from the
surgery presentation theorem, without any appeal to open books.
In Section~\ref{section:HF} I offer the reader an {\em amuse gueule\/}
illustrating the use of contact surgery in Heegaard Floer theory. 
Surgery diagrams also play a supporting role in a contact topological
argument for computing the diffeotopy group of the $3$-manifold
$S^1\times S^2$, as will be explained in Section~\ref{section:diffeotopy}.
An example how contact surgery can be used to detect
so-called non-loose (or exceptional) Legendrian knots
will be given in Section~\ref{section:non-loose}.
Finally, in Section~\ref{section:five} I allow the reader a glimpse of some
recent results in collaboration with Fan Ding and Otto van Koert
on the diagrammatic representation of $5$-dimensional contact
manifolds. The Giroux correspondence reduces the description
of such manifolds to that of a page of an open book
(here: a $4$-dimensional Stein manifold) and the monodromy of the open
book. The $4$-dimensional Stein manifold, in turn, can be described
by a surgery picture that describes the attachment of Stein handles;
the attaching circles for the $2$-handles are Legendrian knots, which
can be visualised in terms of their front projection from $S^3$ (with a
point removed) to a $2$-plane.
In conclusion, one obtains an essentially $2$-dimensional
representation of a contact $5$-manifold. This has implications
on the classification of
subcritically Stein fillable contact $5$-manifolds.
\section{Basic notions and results in contact geometry}
\label{section:basic}
Here I want to recall some fundamental concepts of the
subject. I also mention a few classification and structure theorems
necessary for understanding or putting into perspective
the more recent results described in the subsequent sections.
\subsection{Tight vs.\ overtwisted}
\label{section:tight-ot}
We begin with a dichotomy of contact structures that is specific
to dimension three. A smooth knot $L$ in a contact
$3$-manifold $(M,\xi )$ is called {\em Legendrian\/} if it is everywhere
tangent to the contact structure. If $L$ is homologically trivial in~$M$,
one can find an embedded surface $\Sigma\subset M$ with 
boundary $\partial\Sigma =L$, a so-called {\em Seifert surface\/} for~$L$.
Then $L$ has two distinguished framings (i.e.\ trivialisations of its
normal bundle, which can alternatively be described
by a vector field along and transverse to $L$, or by
a parallel curve obtained by pushing $L$ in the direction of that
vector field): the {\em surface framing}, given by 
a vector field tangent to the surface~$\Sigma$, and the
{\em contact framing}, given by a vector field tangent to the contact
structure~$\xi$. (The surface framing turns out to be independent of the
choice of Seifert surface.)

An embedded $2$-disc $\Delta\subset M$ in
a contact $3$-manifold $(M,\xi )$ is called {\em overtwisted},
if the boundary $\partial \Delta$ is a Legendrian curve whose
contact framing coincides with the surface framing. If one wishes,
one may then arrange that $T_x\Delta =\xi_x$ for all $x\in\partial\Delta$.

A contact $3$-manifold is called {\em overtwisted\/} if it
contains an overtwisted disc; otherwise it is called {\em tight}.
It was shown by Eliashberg~\cite{elia89} that the classification of
overtwisted contact structures on closed
$3$-manifolds is a purely homotopical problem:
each homotopy class of tangent $2$-plane fields contains a
unique overtwisted contact structure (up to isotopy).
For a detailed exposition of Eliashberg's proof see
\cite[Chapter~4.7]{geig08}.

\begin{ex}
Let $(z,r,\varphi )$ be cylindrical coordinates on~$\R^3$.
The contact structure $\xi_{\mathrm{ot}}=\ker (\cos r\, \rmd z+
r\sin r\, \rmd\varphi)$ is an overtwisted contact structure;
each disc $\Delta_{z_0}=\{ z=z_0,\, r\leq \pi\}$ is overtwisted.
\end{ex}

The classification of tight contact structures, on the other
hand, is a very intricate problem that has not yet been solved completely.
It was shown by Bennequin~\cite{benn83}, {\em avant la lettre},
that the standard contact structure $\xist =\ker (\rmd z+x\, \rmd y)$ on
$\R^3$ is tight. We shall return to the classification of tight structures
in Section~\ref{section:space}.
\subsection{Symplectic fillings}
\label{section:filling}
A contact manifold $(M^{2n-1},\xi=\ker\alpha)$ with a
{\em cooriented\/} contact structure is naturally oriented by the
volume form $\alpha\wedge (\rmd\alpha )^{n-1}$. Likewise,
a symplectic manifold $(W^{2n},\omega )$, i.e.\ with $\omega$ a closed
non-degenerate $2$-form, is naturally oriented by the
volume form~$\omega^n$.

\begin{defn}
(a) A compact symplectic manifold $(W^{2n},\omega )$
is called a {\em weak (symplectic) filling\/} of
$(M^{2n-1},\xi=\ker\alpha )$ if $\partial W=M$ as oriented manifolds
and $\omega^{n-1}|_{\xi}> 0$. Here $\partial W$ is oriented by the `outward
normal first' rule.

(b) A compact symplectic manifold $(W^{2n},\omega )$
is called a {\em strong (symplectic) filling\/}
of $(M^{2n-1},\xi=\ker\alpha )$ if $\partial W=M$ and there is a
Liouville vector
field $Y$ defined near $\partial W$, pointing outwards
along $\partial W$, and satisfying $\xi =\ker (i_Y\omega |_{TM})$
(as cooriented contact structure). In this case we say that
$(M,\xi )$ is the {\em convex\/} boundary of $(W,\omega )$.
Here {\em Liouville vector field\/} means that the Lie derivative
$\calL_Y\omega$ --- which is the same as $\rmd (i_Y\omega )$ because of
$\rmd\omega =0$ and Cartan's formula
$\calL_Y=i_Y\circ \rmd +\rmd\circ i_Y$ --- equals~$\omega$.

(c) A {\em Stein filling\/} $(W,J)$ of $(M,\xi )$ is a sublevel set
of an exhausting strictly plurisubharmonic function on a
Stein manifold such that $M=\partial W$ is the
corresponding level set and $\xi$ coincides with the complex
tangencies $TM\cap J(TM)$.
\end{defn}

For the details of (c) I refer the reader to~\cite[Chapter~5.4]{geig08}.
The following implications hold for contact structures:
\[ \mbox{\rm Stein fillable}\;\Longrightarrow\;
\mbox{\rm strongly fillable}\;\Longrightarrow\;
\mbox{\rm weakly fillable}\;\Longrightarrow\;
\mbox{\rm tight}
\]
The first implication is fairly straightforward; the second one is obvious.
For the third implication and references to examples that the converse
implications fail, in general, see~\cite[Chapter~5]{geig08}.
\subsection{Topology of the space of contact structures}
\label{section:space}
Let $M$ be a closed (i.e.\ compact without boundary) 
odd-dimensional manifold.
The space $\Xi (M)$ of contact structures on
$M$ is an open (possibly empty) subset of the space of all
differential $1$-forms on~$M$. According to the Gray
stability theorem~\cite[Theorem~2.2.2]{geig08}, any
smooth homotopy of contact structures on $M$ is induced by
an isotopy of the manifold. So the isotopy classification
of contact structures on $M$ amounts to determining the
set $\pi_0(\Xi (M))$ of path-components.

For $\dim M=3$ it is opportune, thanks to Eliashberg's classification of
overtwisted contact structures, to restrict attention
to tight contact structures. Moreover, we observe that
the sign of $\alpha\wedge \rmd\alpha$ is independent of the choice of
contact form $\alpha$ defining a given contact structure~$\xi$.
If $M$ is oriented, we can thus speak of {\em positive\/} and
{\em negative\/} contact structures. In what follows a choice of
orientation for $M$ will be understood (or specified, if $M$ does
not admit an orientation-reversing diffeomorphism), and we only consider
positive contact structures.

Here are examples of $3$-manifolds with a unique tight contact structure
(up to isotopy); we call this structure the standard contact structure on the
respective manifold and denote it by~$\xist$:
\begin{itemize}
\item $S^3\subset\R^4$, $\xist=\ker (x\, \rmd y-y\, \rmd x+z\, \rmd
t-t\, \rmd z)$,
\item $S^1\times S^2\subset S^1\times\R^3$,
$\xist =\ker (z\, \rmd\theta +x\, \rmd y-y\, \rmd x)$, and
\item $\R^3$, $\xist =\ker (\rmd z+x\, \rmd y)$.
\end{itemize}

\begin{rem}
The standard contact structure on $S^3$, when restricted to
the complement of a point, equals the standard contact structure
on~$\R^3$~\cite[Proposition 2.1.8]{geig08}.
\end{rem}

These results
are due to Eliashberg, cf.~\cite[Chapter~4.10]{geig08}.
Etnyre and Honda~\cite{etho01} have shown that the Poincar\'e homology
sphere $P$ with the opposite of its natural orientation does not admit a
tight contact structure. From a splitting theorem for tight
contact structures due to Colin~\cite{coli97}, cf.~\cite{dige07},
it follows that the connected sum of two copies of $P$,
one with its natural and one with the opposite orientation, does not
admit any tight contact structure for either orientation.

On the $3$-torus $T^3$ the contact structures
$\xi_n=\ker (\sin (n\theta )\, \rmd x+\cos (n\theta)\, \rmd y)$,
$n\in\N$, constitute a complete list
(without repetition) of the tight contact structures {\em up to
diffeomorphism}. The classification up to isotopy is a little more
subtle, see~\cite{elpo94}. As tangent $2$-plane fields, however,
the $\xi_n$ are all homotopic to $\ker \rmd\theta$. This is an instance of
a general phenomenon for toroidal manifolds, i.e.\ manifolds admitting an
embedding of a $2$-torus that induces an injection on fundamental groups:
all such manifolds admit infinitely many tight contact structures.

On the other hand, there are the following finiteness results,
due to Colin--Honda--Giroux~\cite{chg09}:

\begin{itemize}
\item[-] On each closed, oriented $3$-manifold there
are only finitely many homotopy classes of tangent $2$-plane fields
that contain a tight contact structure.
\item[-] Unless the $3$-manifold is toroidal, there are only finitely many
tight contact structures up to isotopy.
\end{itemize}

For $\dim M\geq 5$ there are some existence results for contact structures,
cf.~\cite{geig01a}, but no complete classification on any
contact manifold. An interesting result in this
context is due to Seidel~\cite[Corollary~6.8]{seid07}:
the isomorphism problem for
simply connected closed contact manifolds is
algorithmically unsolvable --- beware, though, that this does not rule
out the practical solution of the problem for a given manifold.

A few things are known about the fundamental group $\pi_1(\Xi (M))$
with a chosen basepoint. For instance, for each $n\in\N$
the group $\pi_1(\Xi(T^3),\xi_n)$
contains an infinite cyclic subgroup~\cite{gego04}. Or,
as shown in~\cite{dige}, the component of $\Xi (S^1\times S^2)$
containing the unique tight contact structure $\xist$ has fundamental
group isomorphic to~$\Z$.
For results about higher homotopy groups of $\Xi (M)$
for higher-dimensional contact manifolds $M$ see~\cite{bour06}.

These results are intimately connected with the topology of
the group $\Diff_0 (M)$ of diffeomorphisms of $M$ isotopic to the
identity. Write $\Xi_0(M)$ for the component of $\Xi (M)$
containing a chosen contact structure~$\xi_0$. Then the map
\[\begin{array}{ccc}
\Diff_0(M) & \longrightarrow & \Xi_0(M)\\
f          & \longmapsto      & Tf(\xi_0)
\end{array}\]
is a Serre fibration, the homotopy lifting property being a
consequence of Gray stability. The fibre $\Cont_0 (M)$,
which need not be connected, consists of those contactomorphisms
(i.e.\ diffeomorphisms that preserve~$\xi_0$) that
are are isotopic (as diffeomorphisms)
to the identity. Thus, the homotopy exact sequence of this Serre
fibration allows us to translate homotopical information about
two of the three spaces $\Cont_0$, $\Diff_0$, and $\Xi_0$
into information about the third.

For the mentioned result $\pi_1(\Xi (S^1\times S^2),\xist )\cong\Z$,
the homotopy type of the topological group $\Diff_0(S^1\times S^2)$
is taken as a given. But there are also examples
where contact topology can be used to extract
information about the diffeomorphism group, see
Section~\ref{section:diffeotopy} below.
\subsection{Convex hypersurfaces}
\label{section:convex}
The notion of a convex hypersurface has been introduced into contact geometry
by Giroux~\cite{giro91}.

\begin{defn}
A vector field $X$ on a contact manifold $(M,\xi )$
is called a {\em contact vector field\/} if its flow preserves
the contact structure~$\xi$. When $\xi$ is written as
$\xi =\ker\alpha$, the condition on $X$ can be stated as $\calL_X\alpha
=\mu\alpha$ for some smooth function $\mu\co M\rightarrow\R$.

A hypersurface $\Sigma\subset M$ is called {\em convex\/}
if there is a contact vector field defined near and transverse to~$\Sigma$.
\end{defn}

\begin{ex}
On $S^1\times\R^2$ with contact structure $\xi =\ker (\cos\theta\, \rmd x-
\sin\theta\, \rmd y)$, the circle $L=S^1\times\{ 0\}$ is Legendrian,
$X=x\,\partial_x+y\,\partial_y$ is a contact vector field, and
$\Sigma =S^1\times\partial D^2$ is a convex surface.
This is actually the universal model for the neighbourhood of
a Legendrian knot in a contact $3$-manifold.
\end{ex}

Convex hypersurfaces, notably in $3$-dimensional contact
manifolds, play an important role in the classification of
contact structures and topological constructions such as surgery.
The reason is the following.

Given a surface $\Sigma$ in a contact
$3$-manifold $(M,\xi)$, the intersection $T\Sigma\cap\xi$
defines a singular $1$-dimensional foliation $\Sigma_{\xi}$
on $\Sigma$, the so-called {\em characteristic foliation}.
Singularities occur at points $x\in\Sigma$ where the tangent plane
$T_x\Sigma$ coincides with the contact plane~$\xi_x$.
It can be shown that the characteristic foliation $\Sigma_{\xi}$
determines the germ of $\xi$ near $\Sigma$. This permits,
for instance, the gluing of contact manifolds along surfaces with the
same characteristic foliation.

In general, the characteristic foliation is difficult to
control. For convex surfaces, however, it turns out that
all the essential information is contained in the
{\em dividing set}, which is defined as the set of points
in $\Sigma$ where the contact vector field is contained
in the contact plane; in a closed surface
this set is a collection of embedded circles. The characteristic foliations
of two convex surfaces with the same dividing set can be made to
coincide after a $C^0$-small perturbation.
\subsection{Open book decompositions}
\label{section:open-book}
Given a topological space $W$ and a homeomorphism $\phi\co W
\rightarrow W$, the
{\bf mapping torus}
$W(\phi )$ is the quotient space
obtained from $W\times [0,2\pi ]$ by identifying $(x,2\pi )$ with
$(\phi (x),0)$ for each $x\in W$. If $W$ is a differential manifold and
$\phi$ a diffeomorphism equal to the identity
near the boundary $\partial W$, then
$W(\phi )$ is in a natural
way a differential manifold with boundary $\partial W\times S^1$.

\begin{figure}[h]
\labellist
\small\hair 2pt
\pinlabel $S^1$ [r] at 237 105
\pinlabel $B$ [l] at 378 14
\pinlabel $\frakp^{-1}(\varphi)$ [l] at 386 435
\endlabellist
\centering
\includegraphics[scale=0.4]{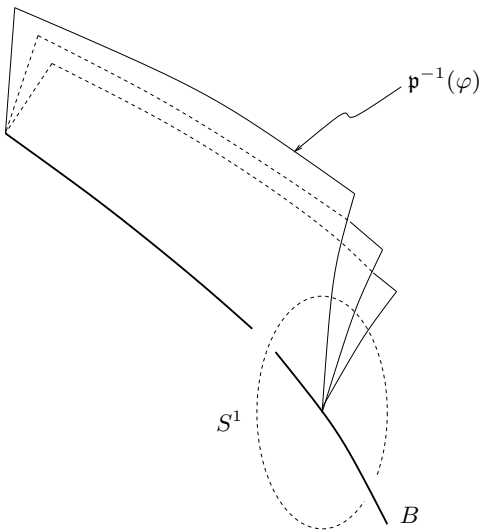}
  \caption{An open book near the binding.}
  \label{figure:book}
\end{figure}

According to an old theorem of Alexander,
cf.~\cite{etny06},
any closed, connected, orientable $3$-manifold can be written in the
form
\[ M_{(\Sigma ,\phi )}:=
\Sigma (\phi )\cup_{\mbox{\rm\scriptsize id}}(\partial\Sigma\times D^2),\]
with $\Sigma$ a compact, orientable surface with boundary;
it can be arranged that the boundary $\partial\Sigma$ is connected
(i.e.\ a single copy of~$S^1$).
Write $B\subset M$ for the link (i.e.\
collection of knots) corresponding to $\partial\Sigma\times\{ 0\}$
under this identification. Then we can define a smooth, locally
trivial fibration
$\frakp\co M\setminus B\rightarrow S^1=\R / 2\pi\Z$ by
\[ \frakp ([x,\varphi ])=[\varphi ]\;\;\;\mbox{\rm for}\;\;\; [x,
\varphi]\in\Sigma (\phi )\]
and
\[ \frakp (\theta ,r\rme^{\bfi\varphi})=[\varphi ]\;\;\;\mbox{\rm for}\;\;\;
(\theta ,r\rme^{\bfi\varphi})\in\partial\Sigma\times D^2\subset\partial
\Sigma\times \C .\]
In other words $B\subset M$ has a tubular neighbourhood of the form
$B\times D^2$, where the fibration $\frakp$ is given by the projection
onto the angular coordinate in the $D^2$-factor. Such a fibration is
called an {\em open book decomposition\/} with {\em binding\/} $B$ and
{\em pages\/} the closures of the fibres $\frakp^{-1}(\varphi )$. Notice that 
each page is a codimension~$1$ submanifold of $M$ with
boundary~$B$, see Figure~\ref{figure:book}.

A submanifold $B\subset M$ that arises as the binding of an open book
decomposition is called a {\em fibred link}.

Conversely, from an open book decomposition of $M$ one can derive a
description of $M$ in the form $M_{(\Sigma ,\phi )}$, so we may think of
an open book decomposition as a pair $(\Sigma ,\phi )$.
The diffeomorphism $\phi$
is called the {\em monodromy\/} of the open book.

In the following definition we call a contact (resp.\ symplectic) form
on an oriented manifold {\em positive\/} if the volume form it defines on
the manifold gives the positive orientation.

\begin{defn}
Let $M$ be a manifold with an open book decomposition $(B,\frakp )$,
where $M$ and $B$ are oriented. The pages of the open book
are oriented consistently with their boundary~$B$.
A contact structure $\xi=\ker\alpha$ on $M$ defined by a positive
contact form is said to be
{\em supported\/} by the open book decomposition $(B,\frakp )$ if
\begin{itemize}
\item[(i)] the $2$-form $\rmd\alpha$ induces a positive symplectic form on
each fibre of~$\frakp$, and
\item[(ii)] the $1$-form $\alpha$ induces a positive contact form on~$B$.
\end{itemize}
\end{defn}

Condition (i) is equivalent to the Reeb vector field $R$ of $\alpha$
being positively transverse to the fibres of~$\frakp$. Recall
that $R$ is defined by the conditions $i_R\rmd\alpha\equiv 0$ and
$\alpha (R)\equiv 1$.

\begin{exs}
(1) The standard contact form on $S^3\subset\C^2$ can be written in
polar coordinates as $\alpha =r_1^2\,\rmd\varphi_1+r_2^2\,\rmd\varphi_2$.
Set $B=\{ r_1=0\}$. Then $\frakp\co S^3\setminus B\rightarrow S^1$,
$\frakp (r_1\rme^{\bfi\varphi_1},r_2\rme^{\bfi\varphi_2})=\varphi_1$
defines an open book whose pages are $2$-discs, and whose
monodromy is the identity map. This open book supports
$\xist =\ker\alpha$, since
$\alpha$ restricts to $\rmd\varphi_2$ along the binding~$B$,
and $\rmd\alpha$ to $r_2\, \rmd r_2\wedge\rmd\varphi_2$ on the
tangent spaces to the pages.

($1^+$) Set $B_+=\{ r_1r_2=0\}$, which is the Hopf link in~$S^3$. Then
$\frakp_+\co S^3\setminus B_+\rightarrow S^1$,
$\frakp_+ (r_1\rme^{\bfi\varphi_1},r_2\rme^{\bfi\varphi_2})=\varphi_1
+\varphi_2$ is an open book whose pages are annuli, and whose
monodromy is a right-handed Dehn twist along the core circle
of the annulus. When oriented as the boundary of a single
page, the binding is a {\em positive\/} Hopf link; the annulus is
called a {\em positive Hopf band}.
For details of these claims, and the fact that
this open book also supports~$\xist$,
see~\cite[Example~4.4.8]{geig08}. Notice that the linking number
of an oriented core circle of the annulus with a push-off along that annulus
equals~$-1$.

(2) The $2$-sphere $S^2$ admits an open book decomposition
where the binding consist of the north and the south pole,
the pages are half great circles between the poles, and the
monodromy is the identity. When we cross this picture with $S^1$
we obtain an open book for $S^1\times S^2$ with binding consisting
of two circles, pages equal to annuli, and monodromy equal to
the identity. The standard contact form $z\, \rmd\theta +
x\,\rmd y-y\,\rmd x$ restricts to $\pm\rmd\theta$ along the binding,
and the Reeb vector field $z\,\partial_{\theta}+x\,\partial_y-y\,\partial_x$
is transverse to the interior of the pages. So this open book
supports the standard contact structure~$\xist$.
\end{exs}

It was shown by Thurston and Winkelnkemper~\cite{thwi75} that
any open book decomposition of a $3$-manifold supports
a contact structure. Giroux~\cite{giro02} observed that
the construction carries over to higher dimensions,
provided the page admits an exact symplectic form
$\omega =\rmd\beta$ which makes it a strong symplectic filling of its
boundary, and the monodromy is symplectic; for
details see~\cite[Chapter~7.3]{geig08}.

Giroux has also shown the converse, which is a much deeper result:

\begin{thm}[Giroux]
\label{thm:Giroux}
Every contact structure on a closed manifold is supported by
an open book decomposition whose fibres are Stein manifolds, and whose
monodromy is a symplectomorphism.\hfill\qed
\end{thm}

Moreover, for $3$-dimensional manifolds he has further refined this
correspondence between contact structures and open books. Given an
open book decomposition of a closed $3$-manifold $M$ with
page $\Sigma$ and monodromy $\phi$, one can form a {\em positive
stabilisation\/} by adding a band to $\Sigma$ along $\partial\Sigma$
and composing $\phi$ with a right-handed Dehn twist
along a simple closed curve running once over the band.
This does not change the underlying $3$-manifold~$M$.
Examples (1) and ($1^+$)
above are an instance of this phenomenon. There is then a
one-to-one correspondence between contact structures on $M$
up to isotopy and open book decompositions of $M$ up
to positive stabilisations and isotopy.

An intrinsic view of this positive stabilisation is to say 
that the page $\Sigma$ is replaced by the plumbing of $\Sigma$
with a positive Hopf band; the plumbing is done in a neighbourhood
of a proper arc in $\Sigma$ and in the Hopf band, respectively;
see~\cite{hare82}.

Analogously, there is a {\em negative stabilisation}, corresponding
to a left-handed Dehn twist or a plumbing with a negative
Hopf band. This will play a role in Corollary~\ref{cor:Harer}.
The corresponding open book of $S^3$ has the negative Hopf link
$B_-$ as binding (which equals $B_+$ as a point set, but one of the two
link components gets the reverse orientation), and the open book
book decomposition is given by $\frakp_-\co
S^3\setminus B_-\rightarrow S^1$,
$\frakp_- (r_1\rme^{\bfi\varphi_1},r_2\rme^{\bfi\varphi_2})=\varphi_1
-\varphi_2$.
\section{A surgery presentation of contact $3$-manifolds}
\label{section:surgery}
\subsection{Dehn surgery}
Let $K$ be a homologically trivial knot in a $3$-manifold $M$. Write
$\nu K\cong S^1\times D^2$
for a (closed) tubular neighbourhood of~$K$. On the boundary
$\partial (\nu K)\cong T^2$ of this tubular neighbourhood
there are two distinguished curves:

\begin{itemize}
\item[1.] The meridian $\mu$, defined as a simple closed curve that
bounds a disc in $\nu K$.
\item[2.] The preferred longitude $\lambda$, defined as a simple closed curve
parallel to $K$ corresponding to the surface framing.
\end{itemize}

Given an orientation of $M$, orientations of $\mu$ and $\lambda$ are
chosen such that the tangent direction of $\mu$ followed by the tangent
direction of $\lambda$ at a transverse intersection point of $\mu$ and
$\lambda$ gives the orientation of $T^2$ (as boundary of $\nu K$).

Let $p,q$ be coprime integers. The manifold $M_{p/q}(K)$ obtained from
$M$ by {\em Dehn surgery\/} along $K$ with {\em surgery coefficient\/}
$p/q\in\Q\cup\{\infty\}$ is defined as
\[ M_{p/q}(K):=\overline{M\setminus\nu K}\cup_g S^1\times D^2,\]
where the gluing map $g$ sends the meridian $*\times\partial D^2$
to $p\mu+q\lambda$, i.e.\ a simple closed curve on $T^2$
in the class $p[\mu ]+q[\lambda ]\in H_1(T^2)$. The resulting manifold is
determined up to diffeomorphism by the surgery coefficient (changing
$p,q$ to $-p,-q$ yields the same manifold).

For $p/q=\infty$ the surgery is trivial. If $p/q\in\Z$,
there is a diffeomorphism $S^1\times D^2\rightarrow\nu K$
sending a standard longitude $\lambda_0=S^1\times\{ *\}$
(with some point $*\in\partial D^2$) to~$p\mu +q\lambda$.
This implies that integer Dehn surgery can be described as
cutting out $S^1\times D^2$ and gluing in $D^2\times S^1$ with the
obvious identification of boundaries. If $M$ is thought of
as the boundary of some $4$-manifold~$W$, the surgered manifold
will be the new boundary after attaching a $2$-handle $D^2\times D^2$
to $W$ along~$M$. For that reason, integer Dehn surgery is also called
{\em handle surgery}.
\subsection{Contact Dehn surgery}
Now suppose that $K$ is a Legendrian knot with respect to some
contact structure $\xi$ on $M$. Then we may replace
$\lambda$ by the longitude corresponding to the
contact framing of~$K$. We now consider
Dehn surgery along $K$ with coefficient $p/q$ as before, but we define
the surgery coefficient with respect to the contact framing. Notice that
the two surgery coefficients differ by an integer depending only on
the Legendrian knot~$K$. This integer, the difference between the
contact framing and the surface framing, is called the
{\em Thurston--Bennequin invariant\/} $\tb (K)$ of~$K$.
(Notice that the contact framing is defined for any
Legendrian knot; the surface framing and $\tb$ are only defined for
homologically trivial ones.)

It turns out that for $p\neq 0$ one can always extend the contact
structure $\xi|_{M\setminus\nu K}$ to one on the surgered manifold
in such a way
that the extended contact structure is tight on the glued-in solid torus
$S^1\times D^2$. Moreover, subject to this tightness condition there
are but finitely many choices for such an extension, and for $p/q=1/k$
with $k\in\Z$ the extension is in fact unique. These observations hinge on
the fact that $\partial (\nu K)$ is a convex surface in the
sense of Section~\ref{section:convex}. On solid tori with
convex boundary condition, tight contact structures have been classified
by Giroux~\cite{giro00} and Honda~\cite{hond00}.

We can therefore speak sensibly of {\em contact $(1/k)$-surgery}.
So the contact surgeries that are well defined {\em and\/} correspond to
handle surgeries are precisely the contact $(\pm 1)$-surgeries.

There is also an {\em ad hoc\/} definition for a contact $0$-surgery,
but here the extension over the glued-in solid torus is necessarily
overtwisted, since the contact framing and the surface
framing of a meridional disc coincide.

The notion of contact Dehn surgery was introduced in~\cite{dige04},
and the following surgery presentation of contact
$3$-manifolds is the main result from that paper.

\begin{thm}
\label{thm:dige}
Let $(M,\xi )$ be a closed, connected contact $3$-manifold. Then
$(M,\xi )$ can be obtained from $(S^3,\xist )$ by contact $(\pm 1)$-surgery
along a Legendrian link.
\end{thm}

\begin{proof}[Sketch proof]
According to a theorem of Lickorish and Wallace, $M$ can be
obtained from $S^3$ by surgery along some link. Since the
reverse of a surgery is again a surgery, we may likewise obtain
$S^3$ by surgery along a link in~$M$.

It is possible to isotope that link in $(M,\xi )$ to a Legendrian link.
Then perform the surgeries as contact surgeries. This yields
$S^3$ with some contact structure~$\xi'$.

Now there is an algorithm for turning each contact surgery into
a sequence of contact $(\pm 1)$-surgeries. Moreover, the contact structures
on $S^3$ are known explicitly (the unique tight one~$\xist$,
and an overtwisted one
in each homotopy class of tangent $2$-plane fields). This allows one
to find a further sequence of contact $(\pm 1)$-surgeries that turns
$(S^3,\xi')$ into $(S^3,\xist )$.

In conclusion, we can obtain $(S^3,\xist )$ from $(M,\xi )$ by
contact $(\pm 1)$-surgery along a Legendrian link. The theorem is then
a consequence of the fact that the converse of a contact $(\pm 1)$-surgery
is a contact $(\mp 1)$-surgery. This `cancellation lemma'
is proved as follows, see Figure~\ref{figure:H2}. 

\begin{figure}[h]
\labellist
\small\hair 2pt
\pinlabel $p_1,p_2$ [r] at 243 528
\pinlabel $q_1,q_2$  at 534 239
\pinlabel {identify via} [l] at 525 431
\pinlabel {Liouville flow} [l] at 525 410
\pinlabel {cut out} [l] at 556 296
\pinlabel {glue in} [bl] at 285 504
\endlabellist
\centering
\includegraphics[scale=0.4]{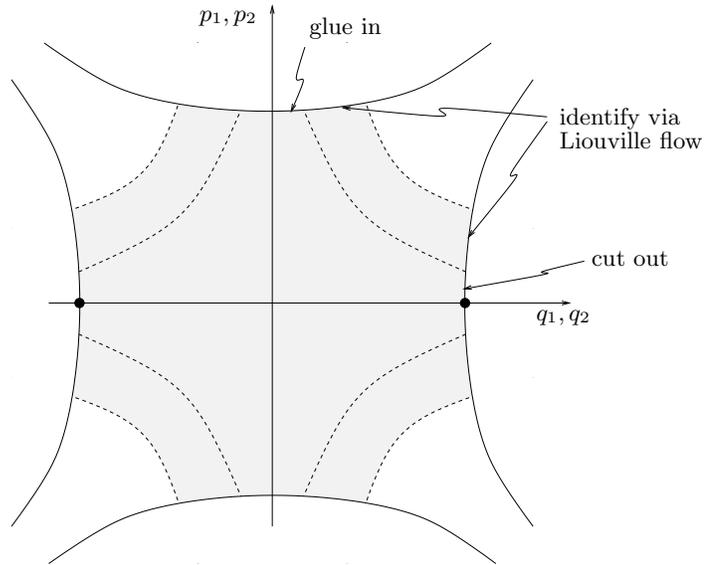}
  \caption{Contact $(-1)$-surgery via Liouville flow.}
  \label{figure:H2}
\end{figure}

Write the Cartesian coordinates on $\R^4$ as $(\bfp ,\bfq )=(p_1,p_2,
q_1,q_2)$. The standard symplectic form on $\R^4$ can then be
written as $\omega =\rmd\bfp\wedge \rmd\bfq :=
\rmd p_1\wedge \rmd q_1+\rmd p_2\wedge \rmd q_2$.
Consider the hypersurfaces $g^{-1}(\pm 1)$, where
$g(\bfp, \bfq )=\bfp^2-\bfq^2/2$, and
the Liouville vector field $Y=2\bfp\, \partial_{\bfp}-\bfq\, \partial_{\bfq}$.
Notice that $Y$ is the gradient vector
field of $g$ with respect to the standard metric on~$\R^4$.
Figure~\ref{figure:H2} gives the local model
for a contact $(-1)$-surgery along the Legendrian circle
$\{ \bfp=0,\; \bfq^2=1\}\subset g^{-1}(-1)$; this follows
from the neighbourhood theorem for Legendrian knots, and
a computation of framings in the local model.

It is clear that the reverse surgery is the one along
$\{ \bfp^2=1,\; \bfq =0\}\subset g^{-1}(1)$ in this
local model, and here a computation of framings shows this
to be a contact $(+1)$-surgery.
\end{proof}

\begin{rem}
In Figure~\ref{figure:H3} of Section~\ref{section:surgery-book}
below we give an alternative description of
contact $(-1)$-surgery that shows how to perform
such a surgery as a symplectic handle
surgery on a (weak or strong) symplectic filling,
so as to obtain a filling of the surgered manifold. This type
of contact surgery had been described earlier by Eliashberg~\cite{elia90}
and Weinstein~\cite{wein91}. Contact $(+1)$-surgery can be interpreted as
a symplectic handle surgery on a {\em concave\/} boundary.
In the `strong' case this means that we have a Liouville
vector field pointing {\em into\/} the filling; in the `weak' case
it is a matter of changing the orientation requirements.
\end{rem}
\section{Applications}
\label{section:applications}

\subsection{From a surgery presentation to an open book}
\label{section:surgery-book}
Given a contact $3$-manifold $(M,\xi )$, Theorem~\ref{thm:dige}
provides us with a Legendrian link $\LL=\LL^-\sqcup\LL^+$
in $(S^3,\xist )$ such that contact $(\pm 1)$-surgery along
the components of $\LL^{\pm}$ yields $(M,\xi )$.
We now want to convert this information into an open book decomposition
of $M$ supporting~$\xi$, which can be done in two steps:

\begin{itemize}
\item[1.] Find an open book for $S^3$ supporting $\xist$, such that
each component of $\LL$ sits on a page of the open book.
\item[2.] Show that contact $(-1)$-surgery (resp.\
$(+1)$-surgery) along
a Legendrian knot sitting on a page of a supporting open book
amounts to changing the monodromy by a right-handed (resp.\
left-handed) Dehn twist.
\end{itemize}

The first step is carried out by Plamenevskaya
in~\cite[Proposition~4]{plam04},
building on work of Akbulut--\"Ozba\u{g}c{\i}~\cite{akoz01}.
The second step is done by Gay~\cite[Proposition~2.8]{gay02}
for contact $(-1)$-surgeries, and for contact surgeries
of both signs by Stipsicz~\cite[pp.~78--79]{stip05}. Their proofs rely
on a result of Torisu about Heegaard splittings of contact
$3$-manifolds along a convex surface into two handlebodies
with a tight contact structure. An alternative proof
of the second step, using only local considerations,
is given by Etnyre~\cite[Theorem~5.7]{etny06};
here I give an independent and self-contained proof.
Like Etnyre's, it is done in a local model,
but the surgery is described
by a smooth model rather than a cut-and-paste procedure.
This proof arose in discussions
with Otto van Koert. Together with Niederkr\"uger he has extended
this argument to higher dimensions; see also~\cite[Proposition~6.2]{akar10}.

In the proposition,
we use the following notation: given a surface $\Sigma$ and
a simple closed curve $L\subset \Sigma$, we write $\rmD_L^+$ for
the diffeomorphism of $\Sigma$ given by a right-handed Dehn twist along~$L$;
a left-handed Dehn twist will be denoted by $\rmD_L^-$.

\begin{prop}
\label{prop:surgery-mono}
Let $(M,\xi =\ker\alpha )$ be a contact $3$-manifold with supporting open book
$(\Sigma ,\phi )$, and let $L$ be a Legendrian knot sitting on a page of
this open book.  Then the contact manifold obtained from $(M,\xi )$ by
contact $(\pm 1)$-surgery along $L$ has a supporting open book
$(\Sigma ,\phi\circ \rmD_L^{\mp})$.
\end{prop}

\begin{proof}
We prove this for a contact $(-1)$-surgery along~$L$; the
case of a contact $(+1)$-surgery is completely analogous.
We begin with a modified local model for a contact $(-1)$-surgery,
see Figure~\ref{figure:H3}. As in Figure~\ref{figure:H2}
we consider $\R^4$ with symplectic form $\omega =
\rmd\bfp\wedge \rmd\bfq$
and Liouville vector field $Y=2\bfp\,\partial_{\bfp}-\bfq\,\partial_{\bfq}$.
But instead of the hypersurface $g^{-1}(-1)$, we now take
the hypersurface $\{\bfq^2=1\}$ as a model for
our contact manifold in a neighbourhood of the Legendrian
knot $L$, which we identify with $\{\bfp =0,\;\bfq^2=1\}$. Perform the surgery
along $L$ by attaching a handle as shown in Figure~\ref{figure:H3},
whose boundary is transverse to the Liouville vector field $Y$
and hence inherits the contact form
$i_Y\omega =2\bfp\, \rmd\bfq +\bfq\, \rmd\bfp$.

\begin{figure}[h]
\labellist
\small\hair 2pt
\pinlabel $q_1,q_2$ [t] at 427 217
\pinlabel $p_1,p_2$ [r] at 217 427
\pinlabel $Y$  at 315 373
\endlabellist
\centering
\includegraphics[scale=0.4]{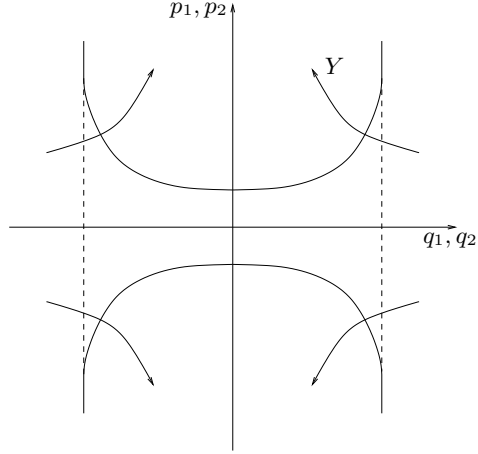}
  \caption{Contact $(-1)$-surgery via handle attachment.}
  \label{figure:H3}
\end{figure}

Consider the map
\[ \begin{array}{rccc}
\frakp\co & \R^4          & \longrightarrow & \R \\
          & (\bfp ,\bfq ) & \longmapsto     & \bfp\bfq.
\end{array} \] 
On the hypersurface $\{ \bfq^2=1\}$ the Reeb vector field $R$
of $i_Y\omega$ takes the form $R=\bfq\, \partial_{\bfp}$, so we
have $R(\frakp )=\bfq^2\equiv 1$ along that hypersurface, which implies
that the Reeb vector field $R$ is transverse to the fibres of the
map~$\frakp$. (These fibres, inside the hypersurface $\{ \bfq^2=1\}$,
are annuli.)
Therefore, by a standard argument involving Gray stability,
cf.~\cite[Chapter~2]{geig08}, we may identify a neighbourhood
of $L\subset M$ with a neighbourhood $\{ \bfp^2<\varepsilon,\; \bfq^2=1\}$
in such a way that $\alpha$ becomes identified with $i_Y\omega$
(restricted to the tangent spaces of the hypersurface $\{\bfq^2=1\}$),
and such that the map $\frakp$ describes
the open book $M\setminus B\rightarrow S^1$ in that neighbourhood.
Notice that the Legendrian knot $L$ lies on the page
$\frakp^{-1}(0)$.


I claim that the map $\frakp$, restricted to the surgered hypersurface
in the local model, still describes an open book 
supporting the contact structure after the surgery. In order to
prove this claim, we need to describe the handle in the model
more explicitly. Following the approach in~\cite{wein91},
we write the surgered manifold in the model as a hypersurface
$\{ F(\bfp^2,\bfq^2)=0\}$, where $F\co\R_0^+\times\R_0^+\rightarrow\R$
is a smooth function with the properties
\[\left\{ \begin{array}{l}
F(0,0)<0,\\[1mm]
\displaystyle{\frac{\partial F}{\partial u}\geq 0,\;\;\;
 \frac{\partial F}{\partial u}>0}\;\;
      \mbox{\rm for}\;\; v=0,\\[3mm]
\displaystyle{\frac{\partial F}{\partial v}\leq 0},\\[3mm]
\displaystyle{\left(\frac{\partial F}{\partial u}\right)^2+
   \left(\frac{\partial F}{\partial v}\right)^2>0},\\[3.5mm]
F(u,1)=0\;\;\mbox{\rm for}\;\; u>\varepsilon^2/4.
\end{array}\right.\]
With $\widetilde{F}(\bfp ,\bfq ):=F(\bfp^2,\bfq^2)$ we have
\[ \rmd\widetilde{F}(Y)=4\bfp^2\frac{\partial F}{\partial u}-
2\bfq^2\frac{\partial F}{\partial v}>0\;\;\mbox{\rm along}\;\;
\{ \widetilde{F}=0\} ,\]
so $\{ \widetilde{F}=0\}$ is indeed a hypersurface transverse to $Y$
that coincides with $\{\bfq^2=1\}$ for $|\bfp| >\varepsilon /2$.

The Reeb vector field $R$ of the contact form induced by $i_Y\omega$
on the hypersurface $\{\widetilde{F}=0\}$ is determined, up to scale,
by the condition that $i_R\rmd (i_Y\omega )=i_R\omega$ be
proportional to~$\rmd\widetilde{F}$. This implies that, up
to a positive factor, the Reeb field is given by
\[ R':=\frac{\partial F}{\partial u}\bfp\,\partial_{\bfq}-
\frac{\partial F}{\partial v}\bfq\,\partial_{\bfp}.\]

From
\[ R'(\frakp )=-\frac{\partial F}{\partial v}\bfq^2+
\frac{\partial F}{\partial u}\bfp^2>0 \;\;\mbox{\rm along}\;\;
\{ \widetilde{F}=0\}\]
it follows that $\frakp$ does indeed define an open book
supporting the contact structure on the surgered manifold.

It remains to verify that this surgery amounts to changing
the monodromy by a right-handed Dehn twist $\rmD^+_L$.
In $3$-manifold topology it is well known that a
Dehn twist on a splitting surface is equivalent to
a surgery along the corresponding curve;
this observation forms the basis for deriving a surgery presentation
of a $3$-manifold from a Heegaard splitting. For completeness I
shall presently provide the argument. {\em A priori}, this only
shows that the surgered manifold admits {\em some\/}
open book decomposition where the monodromy has changed
as described, but not that this is in fact the open book
decomposition given by the map $\frakp$ in the local
model above. A result of Waldhausen~\cite[Lemma~3.5]{wald68}
comes to our rescue: any diffeomorphism of $\Sigma\times [0,1]$
equal to the identity near the boundary is isotopic rel
boundary to a fibre-preserving diffeomorphism; this implies that the
monodromy is determined by a single page and the global topology.
Beware that this is a result specific to dimension~$3$.
Moreover, since I promised a self-contained proof, I show
in Example~(1) following this proof how to give a direct
argument.

\begin{figure}[h]
\labellist
\small\hair 2pt
\pinlabel $\Sigma_0$ [l] at 362 109
\pinlabel $\Sigma_0$ [l] at 362 144
\pinlabel $M_-$ [l] at 362 53
\pinlabel $M_+$ [l] at 362 254
\pinlabel $L$ [t] at 182 139
\endlabellist
\centering
\includegraphics[scale=0.4]{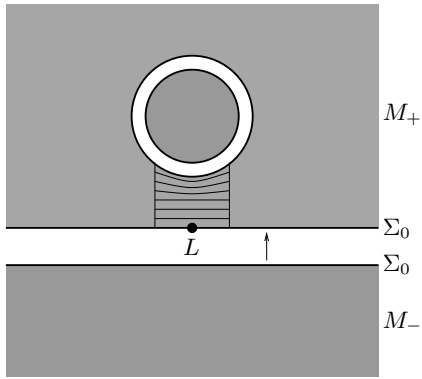}
  \caption{The effect of surgery on the monodromy.}
  \label{figure:surgery-monodromy}
\end{figure}

Imagine that we make an incision in our manifold $M$ along the
page $\Sigma_0$ containing~$L$. Figure~\ref{figure:surgery-monodromy}
shows a cross-section of this incision, orthogonal to~$L$.
In other words, in the figure we see $L$ as a
black dot, and the incision is seen as a horizontal
cut. We think of the positive coorientation to the pages
as pointing up in the figure --- this is the direction of
the flow transverse to the pages that determines the
monodromy. With $M_{\pm}$ we denote neighbourhoods of $L$
on either side of~$\Sigma_0$

If we want to recover the original $M$, we simply reglue
using the identity map. In other words, in our local picture we form
\[ \bigl( M_- +M_+\bigr)/
\bigl( \partial M_-\ni x\sim x\in \partial M_+\bigr) .\]

Now cut $M_+$ open along a $2$-torus lying vertically over~$L$, as
shown in Figure~\ref{figure:surgery-monodromy}. The disc seen in that
figure corresponds to a solid torus~$T$.
Then the right-handed Dehn twist $\rmD_L^+$, which can be thought
of as moving only points in a thin annulus around~$L$,
extends to a diffeomorphism $\rmD_+$ of $M_+\setminus T$ moving only points
in the interior of the
region indicated by (more or less) horizontal lines,
which correspond to annuli,
and acting as a right-handed Dehn twist on each of these annuli.

This diffeomorphism $\rmD_+$ of $M_+\setminus T$, and the identity map
on~$M_-$, induce a diffeomorphism from
\[ \bigl( M_- +(M_+\setminus T)\bigr)/
\bigl( \partial M_-\ni x\sim x\in \partial M_+\bigr) \]
to
\[ \bigl( M_- +(M_+\setminus T)\bigr)/
\bigl( \partial M_-\ni x\sim \rmD_L^+(x)\in \partial M_+\bigr) .\]

So we have changed the monodromy by a right-handed
Dehn twist $\rmD_L^+$, at the price of cutting out a solid torus
and gluing it back after we have performed the diffeomorphism $\rmD_+$
on~$M_+\setminus T$. It remains to show that this cutting and regluing
of $T$ amounts to a $(-1)$-surgery relative to the framing of $L$
given by the page~$\Sigma_0$ (this framing coincides with
the contact framing in the case of a contact structure
supported by the open book).

Think of the meridian $\mu$ on $\partial T$ as the boundary of the
disc seen in Figure~\ref{figure:surgery-monodromy},
oriented clockwise. The longitude $\lambda$ corresponding to
the mentioned framing is a curve parallel to~$L$ (e.g.\ the
curve on~$\partial T$ lying vertically above~$L$).
With the standard orientation of $\R^3$ in our local model,
this longitude points into the picture.

The diffeomorphism $\rmD_-$ has the effect of sending $\lambda$ to
itself and $\mu$ to $\mu +\lambda$. Thus, when we reglue the
solid torus~$T$, its meridian is glued to $\mu-\lambda$.
So this is indeed a $(-1)$-surgery.
\end{proof}

We observed in the above proof that the fibres of $\frakp$ in
the local model are annuli, so it is clear that after the contact
surgery the map $\frakp$ describes an open book
whose monodromy can only have changed by a multiple of
a Dehn twist along the core curve of such a fibre. The same can
be said about the open book obtained from the surgery illustrated in
Figure~\ref{figure:surgery-monodromy}. Therefore, 
the first of the following two examples, where 
the monodromy directly affects the topology, implies
that the change in monodromy is the same in both cases,
i.e.\ a single Dehn twist. This argument allows us to do away
with the reference to Waldhausen.

\begin{exs}
(1) Consider the open book for $S^3$ with binding the positive Hopf link,
with pages diffeomorphic to annuli, and with monodromy a right-handed
Dehn twist along the core circle of the annulus. For any $k\in\Z$
we now want to
find the surgery necessary to turn this into an
open book where the monodromy is a $k$-fold right-handed Dehn twist
(for $k<0$ this means a $|k|$-fold left-handed Dehn twist), i.e.\
we want to add $k-1$ right-handed Dehn twists to the monodromy.
According to the preceding proof, this surgery is given by regluing the
solid torus $T$ by sending its meridian to $\mu -(k-1)\lambda$.
Beware, though, that $\lambda$ does not give the surface framing
(in $S^3$) of the core circle of~$T$.
By Example~(1) in Section~\ref{section:open-book},
the linking number of $L$ with its push-off along the
page is~$-1$. So the surface framing of the
core circle of~$T$ is given by $\lambda'=\mu+\lambda$.
From
\[ \mu -(k-1)\lambda =k\mu - (k-1)\lambda'\]
we deduce that the required surgery is a surgery along an unknot
in $S^3$ with surgery coefficient $-k/(k-1)$. This is the well-known
surgery description of the lens space $L(k,k-1)$,
cf.~\cite[Example~5.3.2]{gost99}.

For an alternative proof that the open book with page an annulus
and monodromy a $k$-fold right-handed Dehn twist is
a lens space $L(k,k-1)$ see~\cite{koni05}. That proof uses Brieskorn
manifolds and generalises to higher dimensions.

(2) The Legendrian unknot
$L=\{ (\rme^{\bfi\varphi},\rme^{-\bfi\varphi})\co
\varphi\in [0,2\pi ]\}$ in $(S^3,\xist )$ is the core circle
in the annulus fibre $\frakp_+^{-1}(0)$ of the open book
$(B_+,\frakp_+)$ supporting~$\xist$, see Example~($1^+$)
of Section~\ref{section:open-book}. As mentioned there, the monodromy
of that open book is a right-handed Dehn twist $\rmD_L^+$ along~$L$.
Thus, when we perform a contact $(+1)$-surgery on~$L$
we obtain the contact structure supported by the open book
with annulus fibres and monodromy equal to
$\rmD_L^+\circ\rmD_L^-=\mathrm{id}$, which by Example~(2)
of Section~\ref{section:open-book} is the standard contact structure on
$S^1\times S^2$.

Notice that $L$ is a standard Legendrian unknot in $(S^3,\xist )$
with Thurston--Bennequin invariant $\tb (L)=-1$ (this characterises
$L$ up to Legendrian isotopy). For alternative proofs that contact
$(+1)$-surgery along $L$ produces $(S^1\times S^2,\xist )$,
see~\cite[Lemma~4.3]{dgs04}, which uses a splitting
along a convex torus, and \cite[Lem\-ma~2.5]{list04}, which uses
the contact invariant from Heegaard Floer theory (see
Section~\ref{section:HF} below). The proof in the present example
is essentially equivalent to that of~\cite[Proposition~4.1]{stip05}.
\end{exs}

The following corollary, in a slightly weaker form, was first proved
by Loi--Piergallini~\cite{lopi01}. In the form presented here,
it is due to Giroux~\cite{giro02}, cf.~\cite[Theorem~5.11]{etny06}.

\begin{cor}[Loi-Piergallini, Giroux]
\label{cor:LPG}
A contact $3$-manifold is Stein fillable if and only
if it admits a supporting open book whose
monodromy is a composition of right-handed Dehn twists.
\end{cor}

\begin{proof}[Sketch proof]
Suppose the contact $3$-manifold $(M,\xi )$ is Stein fillable.
According to a result of Eliashberg, cf.~\cite[Theorem~1.3]{gomp98},
$(M,\xi )$ can be obtained from a connected sum $\# S^1\times S^2$ with its
unique tight contact structure $\xist$ by contact $(-1)$-surgery along a
Legendrian link~$\LL$. There is an open book supporting $\xist$
with trivial monodromy, just as in the preceding example.
One can also construct an open book supporting $\xist$ that
contains $\LL$ on its pages, but may have left-handed Dehn twists
in its monodromy. When we pass to a common stabilisation of these
two open books, we have an open book whose monodromy can be
described by right-handed Dehn twists only {\em and\/} contains
$\LL$ on its pages. Now apply Proposition~\ref{prop:surgery-mono}.

Conversely, suppose that $\xi$ is supported by an open book
$(\Sigma ,\phi )$ with $\phi$ a composition of
right-handed Dehn twists. The contact manifold described
by $(\Sigma ,\mathrm{id})$ admits a Stein filling by the product
$\Sigma\times D^2$; observe that
\[ \partial (\Sigma\times D^2)=(\Sigma\times S^1)\cup
(\partial\Sigma\times D^2)=M_{(\Sigma ,\mathrm{id})} .\]
If the Dehn twists that make up $\phi$ are along homologically
essential curves $L_i$, one can realise each of them as
a Legendrian curve on a page of the open book. By
Proposition~\ref{prop:surgery-mono}, $(M,\xi )$ is then 
Stein fillable as a manifold obtained by
contact $(-1)$-surgery on a Stein fillable manifold.
If an $L_i$ is homologically trivial in~$\Sigma$, one
first writes $\rmD_{L_i}^+$ as a composition of
right-handed Dehn twists along non-separating curves, and then concludes
as before.
\end{proof}

\begin{rem}
There is a related criterion for a contact structure to be tight.
Honda--Kazez--Mati\'c~\cite{hkm07} introduce the notion
of {\em right-veering\/} diffeomorphisms of a surface;
the class of such diffeomorphisms contains those that can be
written as a composition of right-handed Dehn twist. These authors show
that a contact structure is tight if and only if {\em all\/} its
supporting open books have right-veering monodromy.
\end{rem}

The next topological application of contact open books,
due to Giroux--Good\-man~\cite{gigo06}, answers a question of
Harer~\cite[Remark 5.1~(a)]{hare82}.

\begin{cor}[Giroux--Goodmann]
\label{cor:Harer}
Any fibred link in $S^3$ can be obtained from the unknot by finitely many
plumbings and `deplumbings' of Hopf bands.
\end{cor}

\begin{proof}[Sketch proof]
Suppose $B\subset S^3$ is a fibred link, i.e.\ we have an open
book $(B,\frakp )$. We formulate everything in the language of open
books, where the plumbing of a Hopf band corresponds to a positive
or negative stabilisation. Consider the negative
stabilisation $(B_-,\frakp_-)$ of $(B,\frakp )$.

In a negative Hopf band in~$S^3$, the two boundary circles have linking
number $-1$ when oriented as the boundary of the band. Thus, when we
orient the core circle in this band and consider a push-off
of this core circle along the band, with the induced orientation,
their linking number will be~$+1$. It follows that
in the contact structure $\xi_-$ on $S^3$ supported by $(B_-,\frakp_-)$
one can find a Legendrian unknot with $\tb =+1$, which forces
$\xi_-$ to be overtwisted, since such a knot violates
the Bennequin inequality~\cite[Theorem~4.6.36]{geig08}
that holds true in tight contact $3$-manifolds.

Likewise, the unknot in $S^3$ is fibred, and after a negative stabilisation
we obtain an open book $(B',\frakp')$
supporting an overtwisted contact structure.

Once a trivialisation of the tangent bundle of $S^3$ has been chosen,
tangent $2$-plane fields on $S^3$ are in one-to-one correspondence with
maps $S^3\rightarrow S^2$, which are classified by the Hopf
invariant, cf.~\cite[Chapter~4.2]{geig08}.
One can check that a positive stabilisation does
not change the Hopf invariant of the contact structure supported by
the respective open book; the examples (1) and ($1^+$) in
Section~\ref{section:open-book} illustrate this claim.
A negative stabilisation, on the other hand, leads to
a contact structure whose Hopf invariant is one greater.

Thus, by negatively stabilising one of $(B_-,\frakp_-)$ or
$(B',\frakp')$ sufficiently often,
we obtain two open books supporting overtwisted
contact structures $\xi_1,\xi_2$ with the same Hopf invariant.
So $\xi_1,\xi_2$ are
homotopic as tangent $2$-plane fields. By Eliashberg's classification
of overtwisted contact structures, $\xi_1$ and $\xi_2$ are
in fact isotopic as contact structures. Then the Giroux correspondence
guarantees that the underlying open books become
isotopic after a suitable number of further {\em positive\/} stabilisations.
\end{proof}
\subsection{Symplectic caps}
\label{section:caps}
In this section I sketch how Theorem~\ref{thm:dige} can be used
to give a proof of the following theorem, due to
Eliashberg~\cite{elia04} and Etnyre~\cite{etny04}, and then discuss
some of its topological applications. Both Eliashberg and Etnyre base
their proof on an open book decomposition supporting a given
contact structure; the idea for the proof indicated here belongs
to \"Ozba\u{g}c\i\ and Stipsicz~\cite{ozst04}, see~\cite{geig06}
for details.

\begin{thm}[Eliashberg, Etnyre]
\label{thm:caps}
Any weak symplectic filling $(W,\omega )$ of a contact $3$-manifold
$(M,\xi )$ embeds symplectically into a closed symplectic $4$-manifold.
\end{thm}

\begin{proof}[Sketch proof]
We need to show that the given convex filling can be
`capped off', i.e.\ we need to find a concave filling of the contact
$3$-manifold that can be glued to the convex filling so as to produce
a closed $4$-manifold.

The desired cap is constructed in three stages. By Theorem~\ref{thm:dige}
we know that $(M,\xi )$ can be obtained by
performing contact $(\pm 1)$-surgeries on some Legendrian link
$\LL$ in $(S^3,\xist )$. For each component $L_i$ of $\LL$
choose a Legendrian knot $K_i$ in $(S^3,\xist )$ 
with linking number $\lk (K_i,L_i)=1$, and $\lk (K_i,L_j)=0$ for
$i\neq j$. Moreover, we require that $K_i$ have Thurston--Bennequin
invariant $\tb (K_i)=1$, which is the same as saying that
contact $(-1)$-surgery along $K_i$ is the same as a topological
$0$-surgery; such $K_i$ can always be found.
Now attach to $(W,\omega )$ the weak symplectic cobordism $W_1$
between $(M,\xi )$ at the concave end and a new contact manifold
$(\Sigma^3,\xi')$ at the convex end, corresponding
to attaching symplectic handles along the~$K_i$. By our choices,
$\Sigma^3$ will be a homology $3$-sphere.

Thus, after the first step we have embedded $(W,\omega )$ symplectically
into a weak filling $(W\cup_M W_1,\omega')$
of $(\Sigma^3,\xi')$. Since $\Sigma^3$
is a homology $3$-sphere, the symplectic form is exact near
$\Sigma^3=\partial (W\cup_M W_1)$. This allows one to write
down an explicit symplectic form on the cylinder $W_2=\Sigma^3\times [0,1]$
that coincides with $\omega'$ near $\Sigma^3\times\{ 0\}$ and
makes $\Sigma^3\times\{ 1\}$ a {\em strong\/} convex boundary
(with the same induced contact structure~$\xi'$).

We now have a strong filling $(W\cup_M W_1\cup_{\Sigma^3}
(\Sigma^3\times [0,1]),\omega'')$ of $(\Sigma^3,\xi')$.
One could then quote to a result of Gay~\cite{gay02}
that strong fillings can be
capped off; this result, however, is again based on open book
decompositions. Alternatively, we appeal once more to Theorem~\ref{thm:dige},
and argue as follows.
Attach a (strong) symplectic cobordism corresponding to
contact $(-1)$-surgeries that cancel the contact $(+1)$-surgeries
in a surgery presentation of $(\Sigma^3,\xi')$. The new boundary
has a surgery description involving only contact $(-1)$-surgeries
on $(S^3,\xist )$, which implies that it is Stein fillable.
Symplectic caps for Stein fillings have been constructed by
Akbulut--\"Ozba\u{g}c{\i}~\cite{akoz02} and Lisca--Mati\'c~\cite{lima97}.
\end{proof}

This theorem has a number of topological applications, which are
nicely surveyed by Etnyre~\cite{etny09}. For instance,
Kronheimer--Mrowka~\cite{krmr04} used it to show that every
non-trivial knot $K$ in $S^3$ has the (unfortunately named) property~P,
which says that Dehn surgery along $K$ with any surgery coefficient
$p/q\neq\infty$ leads to a $3$-manifold $S^3_{p/q}(K)$ with
non-trivial fundamental group. A more palatable consequence of
this fact is the Gordon--Luecke theorem, which states that knots in
$S^3$ are determined by their complement, cf.~\cite{geig06}.
Theorem~\ref{thm:caps} enters in the Kronheimer--Mrowka proof
as follows. Given a purported counterexample, i.e.\
a non-trivial knot $K\subset S^3$ and some $p/q\neq\infty$
for which $\pi_1(S^3_{p/q}(K))=\{ 1\}$, one constructs with the help
of Theorem~\ref{thm:caps} a certain closed symplectic $4$-manifold that
contains, essentially, $S^3_{p/q}(K)$ as a separating hypersurface.
Deep gauge theoretic results show that such a $4$-manifold cannot
exist.
\subsection{Heegaard Floer theory}
\label{section:HF}
As remarked at the end of Section~\ref{section:surgery},
any contact manifold obtained from a symplectically
fillable contact manifold via contact $(-1)$-surgery will
again be symplectically fillable, and hence in particular tight.
It is not known, in general, whether contact $(-1)$-surgery
on a tight contact $3$-manifold will preserve tightness.
For manifolds with boundary, Honda~\cite{hond02} has an example
where tightness is destroyed by contact $(-1)$-surgery;
for closed manifolds no such example is known.

\begin{figure}[h]
\labellist
\small\hair 2pt
\pinlabel $L$ [tl] at 182 29
\pinlabel $L_{\mathrm{ot}}$ [tl] at 385 48
\pinlabel $+1$ [bl] at 162 133
\endlabellist
\centering
\includegraphics[scale=0.4]{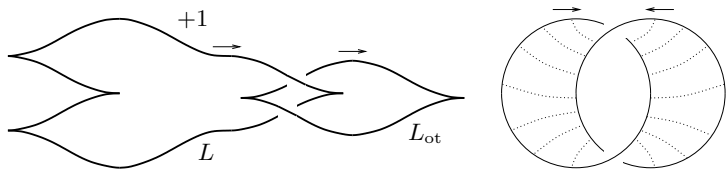}
  \caption{The overtwisted contact manifold $(S^3,\xist )_{+1}(L)$.}
  \label{figure:ot}
\end{figure}

Contact $(+1)$-surgery may well turn
a fillable contact manifold into an overtwisted one.
An example is shown in Figure~\ref{figure:ot} (where the Legendrian knots
in $(\R^3,\ker (\rmd z+x\, \rmd y))\subset (S^3,\xist )$
are represented in terms of their so-called front projection
to the $yz$-plane; the missing $x$-coordinate can be recovered
as the negative slope $x=-\rmd z/\rmd y$). The contact manifold
$(S^3,\xist )_{+1}(L)$ obtained from $(S^3,\xist )$ via contact
$(+1)$-surgery on the `shark' $L$ is overtwisted. Indeed,
the Legendrian knot $L_{\mathrm{ot}}$ bounds an
overtwisted disc in $(S^3,\xist)_{+1}(L)$, as is indicated on the
right side of Figure~\ref{figure:ot}. The
Seifert surface of the Hopf link $L\sqcup L_{\mathrm{ot}}$
shown there glues with a new meridional disc in
in $(S^3,\xist )_{+1}(L)$ to form an embedded disc bounded by
$L_{\mathrm{ot}}$ in the surgered manifold, and the contact framing
of $L_{\mathrm{ot}}$ coincides with the disc framing.

On the other hand, a manifold obtained via contact $(+1)$-surgery
may also be tight, as is shown by example~(2) in
Section~\ref{section:surgery-book}: the tight
contact manifold $(S^1\times S^2,\xist )$
is obtained from $(S^3,\xist )$ by contact $(+1)$-surgery
along a standard Legendrian unknot.

So far the most effective approach towards the question whether
contact $(-1)$-surgery on closed contact $3$-manifolds preserves tightness
comes from the Heegaard Floer theory introduced by
Ozsv\'ath and Szab\'o~\cite{ozsz05}.
Let $(M,\xi )$ be a closed contact $3$-manifold with orientation
induced by the contact structure~$\xi$. We write $-M$ for the manifold
with the opposite orientation. The contact structure $\xi$
determines a natural Spin$^c$ structure $\mathbf{t}_{\xi}$ on~$M$.
Suffice it to say here that Ozsv\'ath and Szab\'o
define a contact invariant $c(M,\xi )$, which lives
in the Heegaard Floer group $\widehat{HF}(-M,\mathbf{t}_{\xi})$,
with the following properties:
\begin{itemize}
\item[-] If $(M,\xi )$ is overtwisted, then $c(M,\xi )=0$.
\item[-] If $(M,\xi )$ is Stein fillable, then $c(M,\xi )\neq 0$.
\end{itemize}

If $(M',\xi')$ is obtained from $(M,\xi)$ by a single contact
$(+1)$-surgery (and hence $(M,\xi)$ by contact $(-1)$-surgery
on $(M',\xi')$), the cobordism $W$ given by the contact
$(+1)$-surgery
induces a homomorphism
\[ F_{-W}\co \widehat{HF}(-M,\mathbf{t}_{\xi})\longrightarrow
\widehat{HF}(-M',\mathbf{t}_{\xi'}).\]
As shown by Lisca and Stipsicz~\cite[Theorem~2.3]{list04}, this homomorphism
maps one contact invariant to the other:
\[ F_{-W}(c(M,\xi ))=c(M',\xi').\]
This immediately implies the following result.

\begin{thm}[Lisca--Stipsicz]
If $c(M',\xi')\neq 0$, then $c(M,\xi )\neq 0$. In particular,
$(M,\xi )$ is tight.\hfill\qed
\end{thm}

In a masterly series of papers, Lisca and Stipsicz have refined this
approach to obtain wide-ranging existence
results for tight contact structures, culminating in their
paper~\cite{list09}, where they give a complete solution to
the existence problem for tight contact structures
on Seifert fibred $3$-manifolds.
\subsection{Diffeotopy groups}
\label{section:diffeotopy}
The {\em diffeotopy group\/} $\D (M)$ of a smooth manifold $M$ is the quotient
of the diffeomorphism group $\Diff (M)$ by its normal subgroup
$\Diff_0(M)$ of diffeomorphisms isotopic to the identity. Alternatively,
one may think of the diffeotopy group as the group $\pi_0(\Diff (M))$
of path components of $\Diff (M)$, since any continuous path
in $\Diff (M)$ can be approximated by a smooth one, i.e.\ an isotopy.

The theorem of Cerf (in its strong form) says that
$\D (S^3)=\Z_2$, that is, up to isotopy there are only
two diffeomorphisms of~$S^3$, the identity and an
orientation reversing one. The diffeotopy groups of a number
of $3$-manifolds are known, for instance those of all lens spaces.

The diffeotopy group $\D (S^1\times S^2)$ was shown to be
isomorphic to $\Z_2\oplus\Z_2\oplus\Z_2$ by Gluck~\cite{gluc62}.
In \cite{dige} we give a contact geometric proof of
Gluck's result. The starting point for this proof is the
uniqueness of the tight contact structure $\xist$ on $S^1\times S^2$.
With Gray stability this easily translates into saying that
any orientation preserving diffeomorphism of $S^1\times S^2$
is isotopic to a contactomorphism of~$\xist$.

In order to find an isotopy of such a contactomorphism $f$ to one
in a certain standard form, and thus to derive Gluck's
theorem, one observes the effect of the contactomorphism
on some Legendrian knot $L$ in $(S^1\times S^2,\xist )$ generating
the homology of $S^1\times S^2$. This can be done in
a contact surgery diagram for $(S^1\times S^2,\xist )$.
The general `Kirby moves' in such a diagram, as
described in~\cite{dige09}, then allow one to find a contact
isotopy from $f(L)$ back to~$L$. This translates into an
isotopy from $f$ to a contactomorphism fixing~$L$. This gives
one enough control over the contactomorphism to
determine its isotopy type.

As an application of such methods, \cite{dige} contains
examples of homologically trivial
Legendrian knots in $(S^1\times S^2,\xist )\#
(S^1\times S^2,\xist )$ that cannot be distinguished by
their classical invariants (i.e.\ the Thurston--Bennequin
invariant and the rotation number, which counts the
rotations of the tangent vector of the Legendrian knot relative
to a trivialisation of the contact structure over a Seifert surface) ---
but which may well be distinguished by performing contact $(-1)$-surgery
on them.
\subsection{Non-loose Legendrian knots}
\label{section:non-loose}
A Legendrian knot $L$ in an overtwisted contact
$3$-manifold $(M,\xi )$ is called {\em non-loose\/} or {\em exceptional\/}
if the restriction of $\xi$ to $M\setminus L$ is tight.
In other words, $L$ has to intersect each overtwisted disc $\Delta$
in $(M,\xi )$ in such a manner that no Legendrian isotopy will
allow one to separate $L$ from~$\Delta$. This is quite a
surprising phenomenon, since overtwisted discs always
appear in infinite families, as in the example given in
Section~\ref{section:tight-ot}.

Exceptional knots were first
described by Dymara~\cite{dyma01}. For a classification of
exceptional unknots in $S^3$ see~\cite[Theorem~4.7]{elfr09}; there
is in fact a unique overtwisted contact structure on $S^3$ that admits
exceptional unknots.

\begin{figure}[h]
\labellist
\small\hair 2pt
\pinlabel $L(n)$ [l] at 453 55
\pinlabel $+1$ [l] at 458 82
\pinlabel $+1$ [l] at 454 107
\pinlabel $-1$ [l] at 453 131
\pinlabel $-1$ [l] at 514 154
\pinlabel $-1$ [l] at 513 201
\pinlabel $n$ [r] at 55 177
\endlabellist
\centering
\includegraphics[scale=0.4]{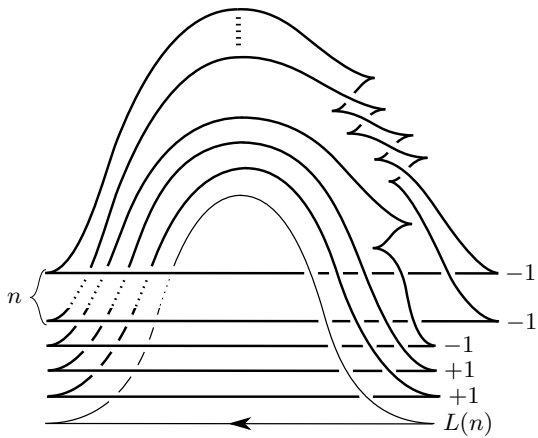}
\caption{Non-loose Legendrian torus knot $T_{2,2n+1}$ in $S^3$
    ($n\geq 1$).}
\label{figure:nonloose}
\end{figure}

Here I want to exhibit an example, due to Lisca
et al.~\cite[Lemma~6.1]{loss09},
which illustrates the use of contact surgery in detecting exceptional
Legendrian knots. Figure~\ref{figure:nonloose} (courtesy of Paolo Lisca
and Andr\'as Stipsicz) shows a surgery
link in $(S^3,\xist )$ (in the front projection); the labels
$\pm 1$ indicate contact $(\pm 1)$-surgeries. The additional
Legendrian knot $L(n)$, which is an unknot in~$S^3$, then represents
a Legendrian knot in the surgered contact manifold~$(M,\xi )$.

By Kirby moves on this surgery diagram one can show that $M$ is
simply another copy of $S^3$, and that $L(n)$ becomes
the torus knot $T_{2,2n+1}$ in this $3$-sphere. Moreover,
with a formula given in \cite[Corollary~3.6]{dgs04},
one can easily compute the Hopf invariant of the contact structure~$\xi$;
it turns out that this differs from the Hopf invariant
$\xist$ on~$S^3$. This implies that $\xi$ and $\xist$ are not homotopic
as tangent $2$-plane fields. Hence, by the uniqueness of the tight
contact structure on~$S^3$, the contact structure $\xi$ must
be overtwisted.

We now want to convince ourselves that $L(n)$ is an exceptional
knot in $(S^3,\xi )$. When we perform contact $(-1)$-surgery
along~$L(n)$, this cancels one of the previous $(+1)$-surgeries.
So the resulting contact manifold is the same as the one
obtained from the original diagram in $(S^3,\xist )$, with
one of the two $(+1)$-surgery knots removed. As seen in
Example~(2) of Section~\ref{section:surgery-book},
a single contact $(+1)$-surgery on a Legendrian unknot
as in Figure~\ref{figure:nonloose} results in $S^1\times S^2$
with its unique tight (and Stein fillable) contact structure.
Further contact $(-1)$-surgeries on this contact manifold
preserve the fillability and hence tightness of the contact structure.
This implies that $T_{2,2n+1}$ is exceptional, for if there were
an overtwisted disc in $S^3\setminus T_{2,2n+1}$, it would survive
to the manifold obtained by surgery along~$T_{2,2n+1}$.

\subsection{Diagrams for contact $5$-manifolds}
\label{section:five}
In the proof of Corollary~\ref{cor:LPG}
we alluded to a result of Eliashberg about the surgery description
of Stein fillable contact $3$-manifolds. That theorem is in fact
a statement about the fillings; in other words, the Stein filling
is obtained by attaching $1$-handles to the $4$-ball (resulting in
a boundary connected sum of copies of $S^1\times D^3$), and then
attaching $2$-handles along Legendrian knots in the boundary
with framing $-1$ relative to the contact framing.

According to Theorem~\ref{thm:Giroux}, any $5$-dimensional
contact manifold $(M,\xi )$ is supported by an open book whose
fibres are Stein surfaces.  By what we just said, those fibres
can be described in terms of a Kirby diagram~\cite{gost99} containing the
information how to attach the $1$- and $2$-handles to the $4$-ball
with its standard Stein structure along its boundary $(S^3,\xist )$.
As in Section~\ref{section:HF}, the pairs of attaching balls
for the $1$-handles and the Legendrian knots along which the
$2$-handles are attached can be drawn in the front projection
of $(\R^3,\ker (\rmd z+x\, \rmd y))$ to the $yz$-plane.

It is not clear how to describe a general symplectic monodromy
in such a diagram. Some monodromies can be
encoded in the diagram, though. For instance, there
are situations where one can `see' Lagrangian spheres in the diagram
(i.e.\ spheres of half the
dimension of the page on whose tangent spaces the symplectic form of the
page vanishes identically),
and one can speak of Dehn twists along such spheres.

Here are some simple examples
with trivial monodromy. Recall from the proof of Corollary~\ref{cor:LPG}
that the manifold $M$ given by an open book with pages $\Sigma$
and trivial monodromy is diffeomorphic to $\partial (\Sigma\times D^2)$.

\begin{exs}
(1) The diagram in Figure~\ref{figure:1-handle} shows a single
$1$-handle; this describes the $4$-manifold $S^1\times D^3$.
So this is a diagram for a contact structure on $\partial
(S^1\times D^3\times D^2)=S^1\times S^4$.

\begin{figure}[h]
\centering
\includegraphics[scale=0.3]{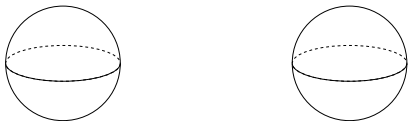}
  \caption{A contact structure on $S^1\times S^4$.}
  \label{figure:1-handle}
\end{figure}

(2) The diagram in Figure~\ref{figure:saucer} shows an
unknot with Thurston--Bennequin invariant $\tb=-1$.
So this corresponds topologically to attaching a $2$-handle with
framing $-2$ relative to the surface framing (given by a spanning disc),
which produces the $D^2$-bundle $\Sigma_{-2}$
over $S^2$ with Euler number~$-2$, see~\cite[Example~4.4.2]{gost99}.
Then $\partial (\Sigma_{-2}\times D^2)$ is the trivial
$S^3$-bundle over~$S^2$. (Observe that the $S^3$-bundles over
$S^2$ are classified by $\pi_1(\SO _3)=\Z_2$; the non-trivial
bundle is detected by the non-vanishing of the second Stiefel--Whitney
class.)

\begin{figure}[h]
\centering
\includegraphics[scale=0.4]{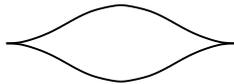}
  \caption{A contact structure on $S^2\times S^3$.}
  \label{figure:saucer}
\end{figure}

(3) In Figure~\ref{figure:shark1} we have an unknot with $\tb =-2$.
So the $4$-manifold encoded by this diagram is the
$D^2$-bundle $\Sigma_{-3}$ over $S^2$ with Euler number~$-3$.
It follows that $\partial (\Sigma_{-3}\times D^2)$ is the unique
non-trivial $S^3$-bundle over~$S^2$, which we write as
$S^2\tilde{\times}S^3$.

\begin{figure}[h]
\centering
\includegraphics[scale=0.4]{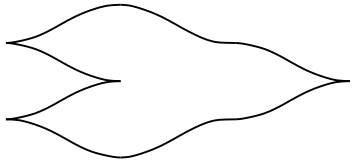}
  \caption{A contact structure on $S^2\tilde{\times}S^3$.}
  \label{figure:shark1}
\end{figure}
\end{exs}

In a forthcoming paper
with Fan Ding and Otto van Koert we exploit the information
contained in such diagrams, and the handle moves
introduced in~\cite{dige}, in order to derive a number of
equivalences of $5$-dimensional contact manifolds. For instance, one
consequence of such moves is that the contact manifold described
by a single Legendrian knot $L$ (and trivial monodromy) will
always be diffeomorphic to $S^2\times S^3$ or $S^2\tilde{\times}S^3$,
and the contact structure is completely determined by the
rotation number of~$L$.

Here is one further observation about open books with trivial
monodromy. From the Seifert--van Kampen theorem one sees that
the fundamental group of 
\[ \partial (\Sigma\times D^2)=\partial\Sigma\times D^2\cup_{\partial}
\Sigma\times S^1\]
is isomorphic to~$\pi_1(\Sigma )$, since any loop in $\partial\Sigma$
is in particular a loop in~$\Sigma$, and $S^1=\partial D^2$
becomes homotopically trivial in~$D^2$. From a Kirby diagram for
$\Sigma$ one can easily read off a presentation of~$\pi_1(\Sigma )$:
each $1$-handle gives a generator, and the attaching circles for
the $2$-handles provide the relations when read as words in
the generators.

As observed by Cieliebak, subcritical Stein fillings (i.e.\
Stein fillings with no handles of maximal index)
split off a $D^2$-factor. Thus, contact manifolds with subcritical Stein
fillings are precisely those admitting an open book with
trivial monodromy.

Combining these two observations, we show in our forthcoming paper
that the contactomorphism type of a subcritically fillable
contact $5$-manifold is, up to a certain stable
equivalence, determined by its fundamental group. This result
is achieved by showing how handle moves in contact surgery
diagrams can be used to effect the so-called Tietze moves
on the corresponding presentation of the fundamental group;
any two finite presentations of a given group are related by such
Tietze moves.

\end{document}